\documentclass[11pt]{article}
\usepackage{amstext,amssymb,amsmath,amsbsy}
\usepackage[dvips,pdfstartview = FitH,colorlinks, urlcolor=green, linkcolor=blue]{hyperref}
\usepackage{amscd}
\usepackage{amsfonts}
\usepackage{indentfirst}
\usepackage{verbatim}
\usepackage{amsmath}
\usepackage{amsthm}
\usepackage{enumerate}
\usepackage{graphicx}
\usepackage{subfig}
\usepackage{color}
\usepackage[OT1]{fontenc}
\usepackage[latin1]{inputenc}
\usepackage[english]{babel}
\usepackage{amssymb}
\usepackage{epsfig}
\newtheorem{theorem}{Theorem}
\newtheorem{lemma}{Lemma}

\newtheorem{remark}{Remark}

\newtheorem{definition}{Definition}
\numberwithin{equation}{section}

\newcommand{\proofend}{\hfill $\Box$ }

\newcommand{\dsp}{\displaystyle}

\newcommand{\eps}{\varepsilon}

\newcommand{\mN}{\mathbb{N}}

\newcommand{\mR}{\mathbb{R}}
\newcommand{\R}{\mathbb{R}}
\newcommand{\mQ}{\mathbb{Q}}

\newcommand{\mZ}{\mathbb{Z}}
\newcommand{\mC}{\mathbb{C}}

\date{\empty}

\title{Dissipative boundary conditions for nonlinear 1-D hyperbolic systems: sharp
conditions through an approach  via time-delay systems}

\author{Jean-Michel Coron\footnote{Sorbonne Universit\'{e}s, UPMC Univ Paris 06, UMR 7598, Laboratoire
Jacques-Louis Lions, F-75005, Paris, France,
E-mail: \texttt{coron@ann.jussieu.fr}. JMC
was supported by ERC advanced grant 266907
(CPDENL) of the 7th Research Framework
Programme (FP7).} \: \:
Hoai-Minh Nguyen \footnote{EPFL SB MATHAA CAMA, Station 8,  CH-1015 Lausanne, Switzerland, E-mail: \texttt{hoai-minh.nguyen@epfl.ch} and School of Mathematics, University of Minnesota, MN, 55455, E-mail: \texttt{hmnguyen@math.umn.edu}. HMN was
supported by NSF grant DMS-1201370, by the Alfred P. Sloan Foundation and by ERC advanced grant 266907
(CPDENL) of the 7th Research Framework
Programme (FP7).}}

\begin{document}
\maketitle
\begin{abstract}
We analyse dissipative boundary conditions for nonlinear
hyperbolic systems in one space dimension. We show  that a previous
known sufficient condition for exponential stability with respect
to the $C^1$-norm is optimal. In particular a
known weaker sufficient condition for exponential stability with respect to the $H^2$-norm is not sufficient
for the exponential stability with respect to the $C^1$-norm. Hence, due to the nonlinearity, even in the case of classical solutions, the exponential stability depends strongly on the
norm considered.
We also give a new sufficient condition
for the exponential stability with respect to the $W^{2,p}$-norm. The methods used are inspired from the
theory of the linear time-delay systems and incorporate the characteristic method.
\end{abstract}

%35L50   	Initial-boundary value problems for first-order hyperbolic systems
%93D20   	Asymptotic stability
\section{Introduction}
Let $n$ be a positive integer. We are concerned with the following nonlinear hyperbolic system:
\begin{equation}\label{sys-P}
u_t + F(u) u_x = 0 \quad  \mbox{ for every } (t, x) \in [0, + \infty) \times [0, 1],
\end{equation}
where $u: [0, + \infty) \times [0, 1] \to \mR^n$ and $F: \mR^n \to {\cal M}_{n, n}(\mR)$. Here, as usual, ${\cal M}_{n, n}(\mR)$ denotes the set of $n \times n$ real matrices. We assume that $F$ is of class $C^\infty$,
$F(0)$ has $n$ distinct real nonzero eigenvalues. Then, replacing, if necessary, $u$ by $Mu$
where $M\in {\cal M}_{n, n}(\mR)$ is a suitable invertible matrix, we may assume that
\begin{equation}\label{cond-F1}
F(0) = \mbox{diag}(\Lambda_1, \cdots, \Lambda_n)
\end{equation}
with
\begin{equation}\label{cond-G1}
\Lambda_i\in \R, \, \Lambda_i \neq \Lambda_j \mbox{ for } i \neq j,
\, i\in\{1,\cdots,n\},\, j\in\{1,\cdots,n\}.
\end{equation}
For simple presentation,  we  assume that,
\begin{equation}\label{lambdai>0}
\Lambda_i > 0 \mbox{ for } i =1, \cdots, n.
\end{equation}
The case where $\Lambda_i$ changes sign can be worked out similarly as in \cite{2008-Coron-Bastin-Novel-SICON}.

\medskip
In this article, we consider the following boundary condition
\begin{equation}\label{bdry-P}
u(t, 0) = G \big( u(t, 1) \big)\quad \mbox{ for every } t \in [0, + \infty),
\end{equation}
where the map $G: \mR^n \to \mR^n$  is of class $C^\infty $ and satisfies
\begin{equation}\label{G(0)=0}
G(0)=0,
\end{equation}
which implies that $0$ is a solution of
\begin{equation}\label{system}
  \left\{
  \begin{array}{ll}
  u_t + F(u) u_x = 0 &\mbox{ for every } (t, x) \in [0, + \infty) \times [0, 1],
  \\
  u(t, 0) = G \big( u(t, 1) \big) &\mbox{ for every } t \in [0, + \infty).
  \end{array}
  \right.
\end{equation}
In this paper, we are concerned about conditions on $G$ for which this equilibrium solution $0$ of \eqref{system} is exponentially stable for \eqref{system}.

\medskip

We first review known results in the linear case, i.e.,  when $F$ and $G$ are linear. In that case,
\eqref{system} is equivalent to
\begin{equation}\label{time-delay-system}
\phi_i(t) = \sum_{j=1}^n K_{ij} \phi_j(t - r_j) \quad \mbox{ for } i =1, \cdots, n,
\end{equation}
 where
 \begin{equation}\label{def-K}
 K =  G'(0) \in {\cal M}_{n \times n}(\mR)
\end{equation}
and
\begin{equation}\label{def-phi-ri}
 \phi_i(t) : = u_i(t, 0), \quad r_i : = 1/ \Lambda_i \quad \mbox{ for } i = 1, \cdots, n.
\end{equation}
Hence \eqref{system} can be viewed as a linear time-delay system. It is known from the work of Hale and Verduyn Lunel \cite[Theorem 3.5 on page 275]{HaleVerduynLunel-book} on delay equations that $0$ is exponentially stable
(in $L^2((0,1);\R^n)$) for \eqref{time-delay-system} if and only if there exists $\delta > 0$ such that
\begin{equation}\label{cond-LN1}
\Big( \mbox{det} \big(Id_n - \big(\mbox{diag} (e^{- r_1 z}, \cdots,  e^{- r_n z})\big)K \big) = 0 , z \in \mC \Big) \implies \Re(z) \le - \delta.
\end{equation}
 For many applications it is interesting to have an exponential stability of \eqref{time-delay-system} which is robust with respect to the small changes on the $\Lambda_i$'s (or, equivalently, on the $r_i$'s),  i.e.,  the speeds of propagation. One says that the exponential stability of $0$ for \eqref{time-delay-system}
 is robust  with respect to the small changes on the $r_i's$ if there exists $\varepsilon \in (0, \text{Min}\{r_1,r_2,\cdots,r_n\})$ such that, for
 every $(\tilde r_1,\tilde r_2,\cdots,\tilde r_n)\in \mR^n$ such that
\begin{equation}\label{Lambda}
  |\tilde r_i- r_i|\leq \varepsilon  \quad \mbox{ for } i =1, \cdots, n,
\end{equation}
$0$ is exponentially stable (in $L^2((0,1);\R^n)$) for
\begin{equation}\label{eqperturb}
  \phi_i(t) = \sum_{j=1}^n K_{ij} \phi_j(t - \tilde r_j) \quad \mbox{ for } i =1, \cdots, n.
\end{equation}
Silkowski (see,  e.g.,  \cite[Theorem 6.1 on page 286]{HaleVerduynLunel-book})
proved that $0$ is exponentially stable (in $L^2((0,1);\R^n)$) for \eqref{time-delay-system} with an exponential stability which is
  robust
with respect to the small changes on the $r_i$'s  if and only if
\begin{equation}\label{cond-LN-2}
\hat \rho_0 \big(K \big) < 1,
\end{equation}
Here
\begin{equation}
\hat \rho_0(K) :  = \max \Big\{ \rho \big( \mbox{diag} (e^{i \theta_1}, \cdots, e^{ i \theta_n}) K \big); \theta_i \in \mR \Big\},
\end{equation}
where, for $M\in {\cal M}_{n \times n}(\mR)$, $\rho(M)$ denotes the spectral radius of $M$. In fact,  Silkowski proved that,
if the $r_i$'s are rationally independent, i.e.,  if
\begin{equation}
\label{rat-ind}
\left(\sum_{i=1}^n q_ir_i=0 \text{ and } q:=(q_1,\cdots, q_n)^T\in \mathbb{Q}^n\right)
\implies \left( q=0\right),
\end{equation}
then $0$ is exponentially stable (in $L^2((0,1);\R^n)$) for \eqref{time-delay-system} if and only if \eqref{cond-LN-2} holds. In \eqref{rat-ind} and in the following, $\mathbb{Q}$
denotes the set of rational numbers.

The nonlinear case has been considered in the literature for more
than three decades. To our knowledge, the first results are due to
Slemrod in \cite{Slemrod} and Greenberg and Li in \cite{GreenbergLi}
in two dimensions, i.e.,   $n=2$. These results were later
generalized for the higher dimensions. All these results rely
on a systematic use of direct estimates of the solutions
and their derivatives along the characteristic curves.
The weakest sufficient condition in this direction was
obtained by  Qin \cite{1985-Qin-Tie-hu},
Zhao \cite{1986-Zhao-Yan-chun} and Li \cite[Theorem 1.3 on page 173]{Li-book}. In these references,
it is proved that $0$  is  exponentially stable for system \eqref{system}
with respect to  the $C^1$-norm if
\begin{equation}\label{cond-Li}
\hat \rho_{\infty} \big(K \big) < 1.
\end{equation}
Here and in the following
\begin{equation}
\hat \rho_p (M): = \inf \big\{ \|\Delta M \Delta^{-1} \|_p; \; \Delta \in {\cal D}_{n, +} \big\} \quad \mbox{ for every } M \in {\cal M}_{n \times n}(\R),
\end{equation}
where ${\cal D}_{n, + }$ denotes the set of all $n \times n$ real diagonal matrices whose entries on the diagonal are strictly positive,  with, for $1 \le p \le \infty$,
\begin{gather}
\label{def|x|p}
\|x \|_p:=\Big(\sum_{i=1}^n|x_i|^p\Big)^{1/p} \quad \forall  x:=(x_1,\cdots,x_n)^T
\in \R^n,\, \forall   p\in [1,+\infty),
\\
\label{def|x|infty}
\|x \|_\infty:=\max\left\{|x_i|;\, i\in \{1,\cdots,n\}\right\} \quad \forall  x:=(x_1,\cdots,x_n)^T
\in \R^n,
\\
\label{def|P|p}
\| M\|_p : = \max_{\|x \|_p = 1} \|M x \|_p\quad \forall M \in {\cal M}_{n \times n}(\R).
\end{gather}
(In fact, in \cite{Li-book,1985-Qin-Tie-hu,1986-Zhao-Yan-chun}, $K$
is assumed to have a special structure; however it is was pointed out in
\cite{2003-De-Halleux-et-al-Automatica} that
the case of a general $K$ can be reduced to the case of this special structure.) We will see later
that \eqref{cond-Li} is also a sufficient condition for the exponential stability with respect
to the $W^{2,  \infty}$-norm (see Theorem~\ref{thm2}). Robustness issues of the exponential stability
 was studied by Prieur et al.
in \cite{2008-Prieur-Winkin-Bastin-MCSS} using again
direct estimates of the solutions
and their derivatives along the characteristic curves.

\medskip
Using a totally different approach, which is based on a Lyapunov stability analysis,
a new criterion on the exponential stability is obtained in \cite{2008-Coron-Bastin-Novel-SICON}: it is proved in this paper that  $0$  is  exponentially stable for system \eqref{system}
with respect to the  $H^2$-norm if
\begin{equation}\label{cond-Coron}
\hat \rho_{2} \big(K\big) < 1.
\end{equation}
This result extends a previous one obtained
in \cite{2007-Coron-Andrea-Novel-Bastin-IEEE} where the same result is established under
the assumption that   $n=2$ and $F$ is diagonal.
See also the prior works \cite{1974-Rauch-Taylor-IUMJ} by Rauch and Taylor, and \cite{02Xu-Sallet} by Xu and Sallet in the case of linear hyperbolic systems.
 It is known (see \cite{2008-Coron-Bastin-Novel-SICON}) that
\begin{equation*}
\hat \rho_0(M) \le  \hat \rho_2 (M) \le \hat \rho_\infty (M)
\end{equation*}
and that the second inequality is strict in general if $n \ge 2$:  for $n \ge 2$ there exists $M \in {\cal M}_{n, n}(\R) $ such that
\begin{equation}\label{compare-1}
\hat \rho_2(M) < \hat \rho_\infty(M).
\end{equation}
In fact, let  $a > 0$ and define
\begin{equation*}
M : = \left( \begin{array}{cc} a & a \\[6pt]
-a & a
\end{array}\right).
\end{equation*}
Then
\begin{equation*}
\hat \rho_2(M) = \sqrt{2} a
\end{equation*}
and
\begin{equation*}
\hat \rho_\infty(M) = 2 a.
\end{equation*}
This implies \eqref{compare-1} in the case $n=2$. The case $n \ge 3$ follows similarly by considering the matrices
\begin{equation*}
\left(\begin{array}{cc}
M & 0 \\[6pt]
0 & 0
\end{array}\right) \in {\cal M}_{n, n}(\R).
\end{equation*}
The Lyapunov approach introduced in \cite{2008-Coron-Bastin-Novel-SICON} has been shown in
\cite{2014-Coron-Bastin} to be applicable to the study the exponential stability with respect to the $C^1$-norm. It gives a new proof that \eqref{cond-Li} implies that $0$  is  exponentially stable for system \eqref{system} with respect to the  $C^1$-norm.

\medskip
The result obtained in \cite{2008-Coron-Bastin-Novel-SICON} is sharp for $n \le 5$. In fact, they
established in \cite{2008-Coron-Bastin-Novel-SICON} the following  result:
\begin{equation*}
\hat \rho_0 = \hat \rho_2 \quad \mbox{ for } n=1,\, 2, \, 3, \, 4, \, 5.
\end{equation*}
For $n \ge 6$, they showed that there exists $M \in {\cal M}_{n, n}(\R)$ such that
\begin{equation*}
\hat \rho_0(M) < \hat \rho_2(M).
\end{equation*}

Taking into account these results, a natural
question is the following: does $\hat \rho_2(K)<1$ implies that $0$  is
 exponentially stable for \eqref{system} with respect to the $C^1$-norm?  We give
a negative answer to this question and prove that
 the condition $\hat \rho_\infty (K)< 1$ is, in some sense, optimal for the  exponential stability with respect to the $C^1$-norm (Theorem~\ref{thm1}).  Hence, different norms require different criteria for the exponential stability with respect to them. Let us emphasize that this phenomenon is due to the nonlinearities: it does not appear when $F$ is constant.  We then show that the condition $\hat \rho_p (K)< 1$ is sufficient to obtain the exponential stability with respect to the $W^{2, p}$-norm (Theorem~\ref{thm2}).  The method used in this paper is strongly inspired from the theory of the linear
 time-delay systems and incorporates the characteristic method.

\medskip
In order to state precisely our first result, we need to recall the compatibility conditions in connection with the
well-posedness for the Cauchy problem associated to \eqref{system}. Let $m\in \mathbb{N}$. Let $\mathcal{H}: C^0([0,1];\R^n)\rightarrow  C^0([0,1];\R^n)$ be a map of class $C^m$. For $k\in \{0,1,\ldots, m\}$, we define,
by induction on $k$, $D^k\mathcal{H}: C^k([0,1];\R^n)\rightarrow  C^0([0,1];\R^n)$ by
\begin{gather}
\label{defDO}
(D^0\mathcal{H})(u):=\mathcal{H}(u) \quad \forall u \in C^0([0,1];\R^n),
\\
(D^k\mathcal{H})(u):=\big((D^{k-1}\mathcal{H}')(u)\big)F(u)u_x  \quad \forall \; u\in  C^k([0,1];\R^n), \, \forall k\in \{0,1,\ldots, m\}.
\end{gather}
 For example, if $m=2$,
\begin{equation}
(D^1\mathcal{H})(u)=\mathcal{H}'(u)F(u)u_x \quad \forall  u\in  C^1([0,1];\R^n),
\end{equation}
\begin{multline}
(D^2\mathcal{H})(u)=\mathcal{H}''(u)\big(F(u)u_x,F(u)u_x\big)+\mathcal{H}'(u)\big(F'(u)F(u)u_x\big)u_x,
\\
+\mathcal{H}'(u)F(u)\big((F'(u)u_x)u_x+F(u)u_{xx}\big)\quad \forall u\in  C^2([0,1];\R^n).
\end{multline}
Let $\mathcal{I}$ be the identity map from $C^0([0,1];\R^n)$ into $C^0([0,1];\R^n)$ and let $\mathcal{G}: C^0([0,1];\R^n)\rightarrow  C^0([0,1];\R^n)$ be defined by
\begin{equation}\label{defcalG}
  \big(\mathcal{G}(v)\big)(x)=G\big(v(x)\big) \quad \text{for every } v \in C^0([0,1];\R^n) \text{ and for every } x\in [0,1].
\end{equation}
 Let $u^0\in C^m([0,1];\R^n)$. We say that
$u^0$ satisfies the compatibility conditions of order $m$ if
\begin{equation}\label{condition-order-m}
  ((D^k\mathcal{I})(u^0))(0)= ((D^k\mathcal{G})(u^0))(1)\quad \text{for every } k\in \{0,1,\ldots, m\} .
\end{equation}
For example, for $m=1$, $u^0\in C^1([0,1];\mR^n)$ satisfies the compatibility conditions of order 1 if and only if
\begin{gather}
\label{compatibilty-C1-0}
u^0(0)=G\big(u(1)\big),
\\
\label{compatibilty-C1-1}
F\big(u^0(0)\big)   u^0_x(0) =  G' \big(u(1) \big) F\big(u^0(1) \big)u^0_x(1).
\end{gather}
With this definition of the compatibility conditions of order $m$, we can recall the following classical theorem
due to Li and Yu \cite[Chapter 4]{LiYu} on the well-posedness of the Cauchy problem associated to \eqref{system}.
\begin{theorem}
\label{wellposedCm}
Let $m\in \mathbb{N}\setminus\{0\}$. Let $T>0$. There exist $\eps>0$ and $C>0$ such that, for every $u^0\in C^m([0,1];\mR^n)$ satisfying the compatibility conditions of order m \eqref{condition-order-m} and such that $\| u^0\|_{C^m([0,1];\R^n)}\leq \eps$, there exists one and only one solution $u\in C^m([0,T]\times [0,1];\R^n)$ of \eqref{system} satisfying the initial condition $u(0,\cdot)=u^0$. Moreover,
%and this solution satisfies
\begin{equation}\label{inedCmCauchy}
  \| u \|_{C^m([0,T]\times [0,1];\R^n)}\leq C \| u^0 \|_{C^m([0,1];\R^n)}.
\end{equation}
\end{theorem}

\begin{remark}
In fact \cite[Chapter 4]{LiYu} is dealing only with the case $m=1$;
however the proof given there  can be adapted to treat the case $m\geq 2$.
\end{remark}

We can now define the notion of exponential stability with respect to the $C^m$-norm.
\begin{definition}
\label{defexpC1}
The equilibrium solution $u \equiv  0$ is exponentially stable for system \eqref{system}  with respect to the $C^m$-norm if there exist $\eps > 0$, $\nu > 0$ and $C>0$ such that, for every $u^0\in C^m([0,1];\mR^n)$ satisfying the compatibility conditions of order m \eqref{condition-order-m} and such that $\| u^0 \|_{C^m([0,1];\R^n)}\leq \eps$, there exists one and only one solution $u\in C^m([0,+\infty)\times [0,1];\R^n)$ of \eqref{system} satisfying the initial condition $u(0,\cdot)=u^0$
and this solution satisfies
\begin{equation*}
\|u(t, \cdot)\|_{C^m([0,1];\R^n)} \le C e^{- \nu t} \| u^0\|_{C^m([0,1];\R^n)} \quad \forall \, t > 0.
\end{equation*}
\end{definition}

With this definition, let us return to the results which are already known concerning the exponential stability with respect to the $C^m$-norm.
\begin{itemize}
\item [(i)] \textbf{For linear $F$ and $G$.} Let $m\in \mathbb{N}$. If $\hat \rho_0\big (G'(0)\big )<1$,  then $0$ is exponentially stable for system \eqref{system}  with respect to the $C^m$-norm and the converse holds if the $r_i$'s are rationally independent. This result was proved for the $L^2$-norm. But the proof can be adapted to treat the case of the $C^m$-norm.
\item [(ii)] \textbf{For general $F$ and $G$.} Let $m\in \mathbb{N}\setminus\{0\}$. If $\hat \rho_\infty \big (G'(0)\big)<1$, then $0$ is exponentially stable for system \eqref{system}  with respect to the $C^m$-norm. This result was proved only for the case $m=1$. However the proofs given in \cite{Li-book, 1985-Qin-Tie-hu, 1986-Zhao-Yan-chun} for this case can be adapted to treat the case $m\geq 2$.
\item [(iii)] \textbf{For general $F$ and $G$, and $n=1$.} Let $m\in \mathbb{N}\setminus\{0\}$. Then $0$ is exponentially stable for system \eqref{system}  with respect to the $C^m$-norm if and only if $\hat \rho_0\big (G'(0)\big )<1$. Note that,
    for $n=1$, the $\hat \rho_p\big (G'(0)\big )$'s
do not depend on $p\in [1,+\infty]$: they are all equal to $|G'(0)|$.
\end{itemize}

The first result of this paper is the following one.
\begin{theorem}\label{thm1} Let $m\in \mathbb{N}\setminus\{0\}$, $n \ge 2$ and $\tau >0$. There exist $F \in C^\infty (\mR^n; {\cal M}_{n \times n}(\mR))$  and  a linear map $G: \mR^n \to \mR^n$  such that $F$ is diagonal, $F(0)$ has distinct positive eigenvalues,
\begin{equation}
\hat \rho_\infty\big( G'(0)\big) < 1 + \tau, \,\hat \rho_0\big( G'(0)\big)=\hat \rho_2\big( G'(0)\big) < 1
\end{equation}
and  $0$ is \textbf{not} exponentially stable for system \eqref{system} with respect to the $C^m$-norm.
\end{theorem}

The second result of this paper is on a sufficient condition for the exponential stability
with respect to  the $W^{2, p}$-norm. In order to state it, we use the following definition, adapted
from Definition~\ref{defexpC1}.

\begin{definition}
\label{defexpW2p} Let $p\in [1,+\infty]$.
The equilibrium solution $u \equiv  0$ is exponentially stable for \eqref{system}  with respect to the $W^{2,p}$-norm if there exist $\eps > 0$, $\nu > 0$ and $C>0$ such that, for every $u^0\in W^{2,p}((0,1);\mR^n)$ satisfying the compatibility conditions of order 1
\eqref{compatibilty-C1-0}-\eqref{compatibilty-C1-1}
and such that
\begin{equation}\label{u0petitC1}
  \|u^0\|_{W^{2,p}((0,1) ;\R^n)} \le \eps,
\end{equation}
there exists one and only one solution $u\in C^1([0,+\infty)\times [0,1];\R^n)$ of \eqref{system} satisfying the initial condition $u(0,\cdot)=u^0$
and this solution satisfies
\begin{equation*}
\|u(t, \cdot) \|_{W^{2,p}((0,1);\R^n)}\le C e^{- \nu t} \| u^0 \|_{W^{2,p}((0,1);\R^n)} \quad \forall \, t > 0.
\end{equation*}

\end{definition}

Again, for every $T>0$,   for every initial condition $u^0\in W^{2,p}((0,1);\mR^n)$ satisfying the compatibility conditions \eqref{compatibilty-C1-0}-\eqref{compatibilty-C1-1} and such that $\| u^0 \|_{W^{2,p}((0,1);\R^n)}$ is small enough, there exist a unique $C^1$ solution $u\in L^\infty([0,T];W^{2,p}((0,1);\R^n))$ of \eqref{system} satisfying the initial condition $u(0,\cdot)=u^0$ (and this solution is in $C^0([0,T];W^{2,p}((0,1);\R^n))$ if $p\in [1,+\infty)$).
The (sketchs of) proof
given in \cite{2008-Coron-Bastin-Novel-SICON} of this result for $p=2$ can be adapted to treat the other cases. Our next result is the following theorem.

\begin{theorem}\label{thm2} Let $p\in [1,+\infty]$. Assume that
\begin{equation}
\hat \rho_{p} \big(G'(0)\big) < 1.
\end{equation}
Then, the equilibrium solution $u \equiv  0$ of the system \eqref{system} is exponentially stable with respect to the $W^{2,p}$-norm.
\end{theorem}
Let us recall that the case $p=2$ is proved in \cite{2008-Coron-Bastin-Novel-SICON}. Let us emphasize that,
even in this case, our proof is completely different from the one given in \cite{2008-Coron-Bastin-Novel-SICON}.

\begin{remark} The notations on various conditions on exponential stability used in this paper are different from the ones in \cite{2008-Coron-Bastin-Novel-SICON}. In fact, one has
\begin{equation*}
\hat \rho_0 = \rho_0,  \quad \hat \rho_2 = \rho_1, \quad \mbox{ and } \quad \hat \rho_\infty = \rho_2.
\end{equation*}
Here $\rho_0$, $\rho_1$, and $\rho_2$ are the notations used in \cite{2008-Coron-Bastin-Novel-SICON}.
\end{remark}

The paper is organized as follows. In Sections~\ref{sect-thm1} and \ref{sect-thm2}, we establish Theorems~\ref{thm1} and \ref{thm2} respectively.

\section{Proof of Theorem~\ref{thm1}}\label{sect-thm1}

We give the proof in the case $n=2$. The general cas $n\geq 2$ follows immediately from the case considered here.

Let $F\in C^{\infty}(\mR^2;\mathcal{M}_{2 \times 2}(\R))$ be such that
\begin{equation}
\label{defF}
F(u) = \left( \begin{array}{cc} \Lambda_1 & 0 \\[6pt]
0 & \displaystyle \frac{1}{ r_2 + u_2}
\end{array} \right) \quad \forall u=(u_1,u_2)^T\in \R^2 \text{ with } u_2>-\frac{r_2}{2},
\end{equation}
for some $0 < \Lambda_1 < \Lambda_2$. We recall that
\begin{equation*}
r_1 = 1/ \Lambda_1 \quad \mbox{ and } \quad r_2 = 1/ \Lambda_2.
\end{equation*}
We assume that $r_1$ and $r_2$ are independent in $\mZ$,  i.e.,
\begin{equation}\label{independence}
\left(k_1 r_1 + k_2 r_2 = 0\text{ and } (k_1,k_2)^T\in \mathbb{Z}^2\right)\implies \left(k_1 = k_2 = 0\right).
\end{equation}
\medskip
Define $G: \mR^2 \to \mR^2 $ as the following linear map
\begin{equation}
\label{defG}
G(u)  :=
a\left(
\begin{array}{cc}
  1   & \xi   \\
  -1   & \eta   \\
\end{array}
\right) u \quad \text{for } u\in \R^2.
\end{equation}
Here $a> 0$ and $\xi, \eta$ are two positive numbers such that
\begin{equation}\label{cond-xieta}
\mbox{ if } P_k(\xi, \eta) = 0  \quad \mbox{ then } \quad P_k \equiv 0,
\end{equation}
for every  polynomial $P_k$ of degree $k$ ($k \ge 0$) with rational coefficients.

Note that if
\begin{equation}
  a \mbox{ is close to } 1/2 \quad \mbox{ and } \quad \xi, \eta \mbox{ are close to 1},
\end{equation}
then
\begin{equation}
\hat \rho_\infty(G) \mbox{ is close to } 1
\end{equation}
and
\begin{equation}
\hat \rho_0(G) = \hat \rho_2(G) \mbox{ are close to } \frac{1}{\sqrt{2}}<1.
\end{equation}
Here, and in the following, for the notational ease, we use the convention  $G = K =  G'(0)$.

\medskip
Let $\tau_0>1$ (which will defined below). We take $a \in \mQ$,  $a>1/2$ but close to $1/2$ and  choose $\xi, \eta >1$ but close to 1 so that
\begin{gather}
\label{rhoinftyGbon}
\hat \rho_\infty (G) < \tau_0,
\\
\label{axietapres}
a(1+\xi+\eta )\leq 2,
\end{gather}
and there exists $ c>0 $  such that
\begin{equation}\label{def-c}
 \frac{\max\{\xi, \eta \}}{a(\xi + \eta)} < c < 1.
\end{equation}
We also impose that $\xi, \, \eta$ satisfy \eqref{cond-xieta}.

We start with the case $m=1$. We argue by contradiction. We assume that there exists $\tau_0 > 1$ such that for all $G$
with $\hat \rho_\infty\big( G'(0) \big) < \tau_0$, there exist $\eps_0$, $C_0$, $\nu$ positive numbers such that
\begin{equation}
\label{udecroitexp}
\| u(t, \cdot ) \|_{C^1([0,1];\R^2)} \le C e^{-\nu t } \|u^0 \|_{C^1([0,1];\R^2)},
\end{equation}
if $u^0\in C^1([0,1];\R^2)$ satisfies the compatibility conditions
\eqref{compatibilty-C1-0}-\eqref{compatibilty-C1-1} and is such that $\| u^0 \|_{C^1([0,1];\R^2)} \le \eps_0$. Here  $u$
denotes the solution of \eqref{system} satisfying
the initial condition $u(0,\cdot)=u^0$.

\medskip

%Note that for such $a$, $\xi, \eta$,
%\begin{equation}\label{G-norm}
%|G |_\infty < 3/2.
%\end{equation}

\medskip
Assume that $u \in C^1([0, + \infty) \times [0, 1];\R^2)$ is a solution to \eqref{system}.
Define
\begin{equation*}
v(t) = u(t, 0).
\end{equation*}
Then
\begin{equation}\label{Delay1}
v\Big(t + r_2 + v_2(t) \Big) = v_1\Big(t + r_2 + v_2(t) - r_1\Big) G_1 + v_2(t) G_2.
\end{equation}
where $G_1$ and $G_2$ are the first and the second column of $G$. Equation \eqref{Delay1} motivates our construction below.

\medskip
Fix $T > 0$ (arbitrarily  large) such that
\begin{equation*}
T - (k r_1 + l r_2) \neq 0\quad  \mbox{ for every } k, l \in \mN.
\end{equation*}
Let $\eps\in (0,1)$ be (arbitrarily) small  such that
\begin{equation}\label{nodefect}
\inf_{k, l \in \mN} |T - (kr_1 + lr_2)| \ge \eps.
\end{equation}
(Note that the smallness of $\eps$ in order to have \eqref{nodefect} depends on $T$: It goes to $0$ as $T\rightarrow +\infty$.) Let $n$ be the integer part of $T/r_2$ plus $1$.  In particular $ n r_2 > T$. Fix
$n$ rational points  $(s_i^0, t_i^0)^T \in \mQ^2$,  $i = 1, \cdots, n$,   such that their  coordinates are distinct, i.e., $s_i^0 \neq s_j^0$, $t_i^0 \neq t_j^0$ for $i \neq j$, and
\begin{equation}\label{initial}
\|(s_i^0, t_i^0)  \|_\infty \le \eps^3/ 4^n \quad  \mbox{ for every } i \in \{1,\cdots,n\}.
\end{equation}
For $0 \le k \le n-1$, we define $(s_i^{k+1}, t_i^{k+1})^T$ for $i=1, n-(k+1)$ by recurrence as follows
\begin{equation}\label{construction-st}
(s_i^{k+1}, t_i^{k+1})^T =  G (s_i^k, t_{i+1}^k)^T = a \left( \begin{array}{c}
s_i^k + \xi t_{i+1}^k  \\[6pt]
- s_i^k + \eta t_{i+1}^k
\end{array}\right).
\end{equation}
%Note that  if $v(t) = (s_i^{k+1}, t_i^{k+1})^T$, $v(t - r_1) = (s_i^{k}, t_i^k)^T$ and $v(t - r_2 - t_{i+1}^k) = (s_{i+1}^k, t_{i+1}^k)$ for some $t$, then $v(t) = v(t-r_1) G_1 + v_2(t- r_2 - t_{i+1}^k) G_2$, which is of form \eqref{Delay1}.

\medskip
Set
\begin{equation}
\label{initVetdV}
V(T) := (s_1^n, t_1^n), \quad dV(T) = \eps (1, 0)^T.
\end{equation}
Define
\begin{gather}
\label{Tinit}
T_{1} := T- r_1, \quad   T_{2} := T - r_2 - t_{2}^{n-1},
\\
\label{Vinit}
V(T_{1}) = (s_1^{n-1}, t_1^{n-1}),  \quad V(T_{2}) = (s_2^{n-1}, t_2^{n-1}),
\\
\label{dVinit}
dV(T_1) =  \eps \Big(\frac{\eta}{a(\xi + \eta)}, 0 \Big),   \quad dV(T_2) = \eps \Big(0, \frac{1}{a(\xi + \eta)} \Big).  \end{gather}

Assume that $T_{\gamma_1  \cdots  \gamma_k}$ is  defined for $\gamma_i = 1, 2$. Set
\begin{equation}\label{def-t1}
T_{\gamma_1  \cdots  \gamma_k  1} = T_{\gamma_1  \cdots  \gamma_k} - r_1
\end{equation}
and
\begin{equation}\label{def-t2}
T_{\gamma_1 \cdots \gamma_k 2} = T_{\gamma_1 \cdots \gamma_k} - r_2 - t_{1 + l}^{n-(k+1)}.
\end{equation}
where \footnote{Roughly speaking, $l$ describes the number of times which comes from $r_2$ in the construction of $\gamma_1 \cdots \gamma_k$. }
\begin{equation}\label{def-l}
l = \sum_{j = 1}^k (\gamma_j - 1).
\end{equation}
Note that, by \eqref{initial}, \eqref{construction-st}, \eqref{Tinit},  \eqref{def-t1}, \eqref{def-t2} and \eqref{def-l}
\begin{equation}\label{estTgamma}
  \Big|T_{\gamma_1\cdots\gamma_k}-kr_1-(r_2-r_1)
  \sum_{j=1}^k\left(\gamma_j-1\right)\Big|\leq C \varepsilon^3 \quad\forall k\in\{1,\cdots,n\},
\end{equation}
for some $C>0$ which is independent of $T>r_1$ and $\varepsilon \in (0,+\infty)$.

\begin{figure}[hbtp]
\begin{center}
\epsfig{file=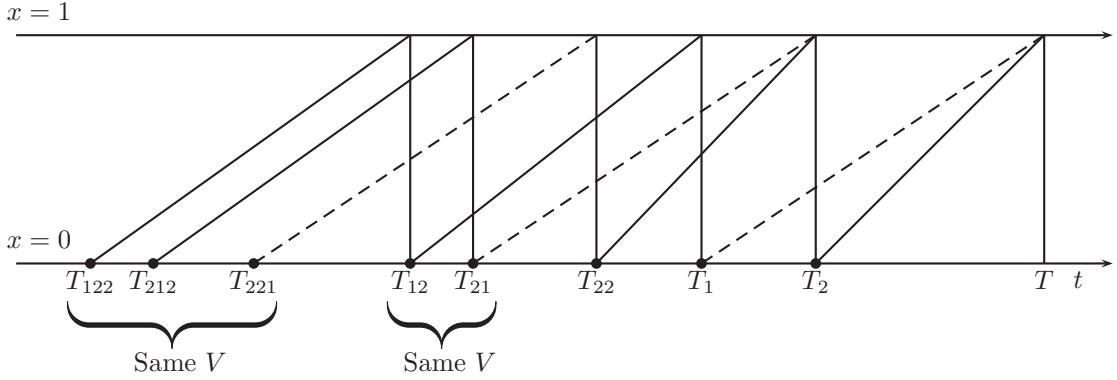, width=\linewidth, clip=}
\caption{$V(T_{122})=V(T_{212})=V(T_{221})\not = V(T_{12})=V(T_{21})$ and the $T_\gamma$'s are different. The slope of the dashed lines is $\Lambda_1=r_1^{-1}$.}
\label{ref-Tgamma-fig}
\end{center}
\end{figure}

We claim that
\begin{equation}\label{T}
\text{the }T_{\gamma_1 \cdots \gamma_k}, \, k\in\{1,\cdots,n-1\}, \mbox{ are  distinct}.
\end{equation}
(See fig. \ref{ref-Tgamma-fig}.) We admit this fact, which will be proved later on, and continue the proof.
%\medskip Therefore, $T_{\alpha_1, \cdots, \alpha_{n-1}}$ are distinct.  This fact plays an important role in our construction.

\medskip
Define $V(T_{\gamma_1  \cdots  \gamma_k \gamma_{k+1}})$ and  $dV (T_{\gamma_1  \cdots  \gamma_k \gamma_{k+1}})$ as follows
\begin{equation}\label{def-v}
V (T_{\gamma_1  \cdots  \gamma_k  \gamma_{k+1}}) = (s_{1 + l}^{n-(k+1)}, t_{1 + l}^{n - (k+1)})^T
\end{equation}
and
\begin{equation}\label{def-dv}
dV (T_{\gamma_1  \cdots  \gamma_k  1}) = (x, 0)^T \quad \quad  dV (T_{\gamma_1  \cdots  \gamma_k  2})
 = (0, y)^T,
\end{equation}
where $l$ is given by \eqref{def-l} and the real numbers $x, y$ are chosen such that
\begin{equation}
\label{defxy}
G (x, y) ^T = dV (T_{\gamma_1  \cdots  \gamma_k}).
\end{equation}
Let us also point that, by \eqref{dVinit} and \eqref{def-dv},
\begin{equation}\label{aumoinsunecompsantenulle}
  \text{at least one of the two components of $dV (T_{\gamma_1  \cdots  \gamma_k})$ is 0.}
\end{equation}
 From \eqref{defG}, we have
\begin{equation}\label{G-1}
  G^{-1}=\frac{1}{a(\eta+\xi)}
  \begin{pmatrix}
  \eta & -\xi
  \\
  1 & 1
  \end{pmatrix}.
\end{equation}
%\begin{equation}\label{observation1}
% G \Big(\frac{\eta}{a (\xi + \eta)}, \frac{1}{a(\xi + \eta)}\Big)^T = (1, 0)^T
%\end{equation}
%and
%\begin{equation}\label{observation2}
% G\Big(-\frac{\xi}{a (\xi + \eta)}, \frac{1}{a(\xi + \eta)}\Big)^T = (0, 1)^T.
%\end{equation}
It follows from  \eqref{def-c}, \eqref{def-dv}, \eqref{defxy}, \eqref{aumoinsunecompsantenulle}
and \eqref{G-1} that
\begin{equation}\label{cond-dv}
\|dV (T_{\gamma_1  \cdots  \gamma_k  \gamma_{k+1}})\|_\infty \le c \| dV (T_{\gamma_1  \cdots  \gamma_k}) \|_\infty.
\end{equation}

\medskip
Using \eqref{T}, we may  construct $\mathfrak{v}\in C^1([0, r_1];\mR^2)$ such that
\begin{equation}\label{def-v-1}
\mathfrak{v}' (T_{\alpha_1 \cdots \alpha_k}) = dV(T_{\alpha_1 \cdots \alpha_k}),
\end{equation}
and
\begin{equation}\label{def-v-2}
\mathfrak{v} (T_{\alpha_1 \cdots \alpha_k}) = V(T_{\alpha_1 \cdots \alpha_k}),
\end{equation}
if
\begin{equation*}
T_{\alpha_1 \cdots \alpha_k} \in (0, r_1),
\end{equation*}
(recall that $r_1 > r_2 > 0$ and $n r_2 > T$).  It follows from \eqref{axietapres}, \eqref{initial}, \eqref{construction-st}, \eqref{def-v} and \eqref{def-v-2} that
\begin{equation}\label{norm-v}
\| \mathfrak{v}(T_{\alpha_1 \cdots \alpha_k})\|_\infty \le \eps^3 \quad \mbox{ if } T_{\alpha_1 \cdots \alpha_k} \in (0, r_1).
\end{equation}
Let $T_{\alpha_1 \cdots \alpha_k} \in (0, r_1)$ and $T_{\gamma_1 \cdots \gamma_m}\in (0, r_1)$ be such that
\begin{equation}
\label{vdifferent}
\mathfrak{v}(T_{\alpha_1 \cdots \alpha_k}) \neq  \mathfrak{v}(T_{\gamma_1 \cdots \gamma_m}).
\end{equation}
 From  \eqref{construction-st}, \eqref{def-v}, \eqref{def-v-2} and \eqref{vdifferent},  we get that
\begin{equation}\label{cnpournonequality}
  k\not = m \text{ or } \mbox{card}\{i\in \{1,\cdots,k\};\, \alpha_i=1\}\not
  = \mbox{card}\{i\in \{1,\cdots,m\};\, \gamma_i=1\}.
\end{equation}
See also Fig. \ref{ref-Tgamma-fig}.

From \eqref{nodefect}, \eqref{Tinit}, \eqref{def-t1}, \eqref{def-t2} and \eqref{cnpournonequality}, we get that, at least if $\varepsilon >0 $ is small enough,
\begin{equation}\label{dist-v}
|T_{\alpha_1 \cdots \alpha_k} - T_{\gamma_1 \cdots \gamma_m}| \ge \eps/ 2.
\end{equation}
Using \eqref{nodefect}, \eqref{norm-v} and \eqref{dist-v}, we may also impose that
\begin{gather}
\label{v=0presde0}
\text{$\mathfrak{v} = 0$ in a neighborhood of $0$ in $[0, r_1]$,}
\\
\label{v=0presder1}
\text{$\mathfrak{v} = 0$ in a neighborhood of $r_1$ in $[0, r_1]$,}
\\
\label{v=0presder2}
\text{$\mathfrak{v} = 0$ in a neighborhood of  $r_2$,}
\\
\label{v-1}
\|\mathfrak{v} \|_{C^1([0, r_1])} \le C \max\{\eps^2, A \},
\end{gather}
where
\begin{equation}
\label{defA}
A : = \max\big\{ \|dV(T_{\alpha_1 \cdots \alpha_k})\|_\infty; \; T_{\alpha_1 \cdots \alpha_k} \in (0, r_1) \big\}.
\end{equation}
In \eqref{v-1},  $C$ denotes a positive constant which does not depend on $T>r_1$ and on $\eps>0$ provided that $\eps>0$ is small enough, this smallness depending on $T$. We use this convention until the end of this section and the constants $C$ may vary from one place to another.

Note that if $T_{\alpha_1 \cdots \alpha_k} \in (0, r_1)$ then
\begin{equation*}
k r_1 > T/2.
\end{equation*}
It follows that
\begin{equation*}
k > T/ (2 r_1),
\end{equation*}
which, together with \eqref{initVetdV}, \eqref{cond-dv} and $c\in(0,1)$, implies that
\begin{equation}
\label{dVpetit}
\|dV(T_{\alpha_1 \cdots \alpha_k}) \|_\infty \le \eps c^{T/ (2 r_1)}.
\end{equation}
 From \eqref{v-1} and \eqref{dVpetit}, one has
\begin{equation}\label{v-2}
\| \mathfrak{v}\|_{C^1([0, r_1];\R^2)} \le C \max\big\{\eps^2,  \eps c^{T/ (2 r_1)} \big\}
\leq C\eps c^{T/ (2 r_1)}.
\end{equation}

Let $\tilde u\in C^1([0,r_1]\times [0,1];\R^2)$ be the solution to the backward Cauchy problem
\begin{equation}\label{tildeubackward}
    \left\{
  \begin{array}{ll}
  \tilde u_t + F(\tilde u) \tilde u_x = 0 &\mbox{ for every } (t, x) \in [0,r_1] \times [0, 1],
  \\
  \tilde u(t, 1) = G^{-1} v(t) &\mbox{ for every } t \in [0,r_1],
  \\
  \tilde u(r_1,x)=0 &\mbox{ for every } x \in [0,1].
  \end{array}
  \right.
\end{equation}
Note that, by \eqref{v=0presder1}, the boundary condition  at $x=1$ for the backward Cauchy problem \eqref{tildeubackward}
vanishes in a neighborhood of $r_1$ in $[0,1]$ and therefore the necessary compatibility conditions for the existence of $\tilde u$, namely
\begin{equation}\label{bondarycompatildeu}
G^{-1} v(t_1)=0 \text{ and } G^{-1} v'(t_1)=0,
\end{equation}
are satisfied. Moreover, if $\varepsilon>0$ is small enough this solutions indeed exists by \cite[pp. 96-107]{LiYu}. Let $u^0 \in C^1([0, 1];\R^2)$ be defined by
\begin{equation}
\label{defu}
u^0(x) := \tilde u(0,x)\quad  \mbox{for every } x \in [0, 1].
\end{equation}
Using \eqref{v-2}, \eqref{def-c} and the definition of $u^0$, we have
\begin{equation}\label{choice-u0}
\| u^0 \|_{C^1([0,1];\R^2)} \le C \| v \|_{C^1([0, r_1];\R^2)} \le C  \max\big\{\eps^2,  \eps c^{T/ (2 r_1)} \big\}
\le C \eps.
\end{equation}
Note that $u^0$ satisfies the the compatibility condition \eqref{compatibilty-C1-0} and \eqref{compatibilty-C1-1} since,
by \eqref{v=0presder1} and \eqref{v=0presder2},
$u^0$ vanishes in a neighborhood of $0$  in $[0,1]$ and, by \eqref{v=0presde0}, $u^0$ vanishes in a neighborhood of $1$  in $[0,1]$. Let $u\in C^1([0,+\infty)\times[0,1];\R^2)$
be the solution of \eqref{system} satisfying the initial condition
\begin{equation*}
u(0,x) = u^0(x) \quad \mbox{ for every } x\in [0,1].
\end{equation*}
Since $0$ is assumed to be exponentially stable for \eqref{system} with respect to the $C^1$-norm,  $u$ exists for all positive time if $\eps$ is small enough. Let us define $v\in C^1([0,+\infty);\R^2)$ by
\begin{equation}\label{deffrakv}
v(t):=u(t,0)\quad \mbox{ for every } t\in[0,+\infty).
\end{equation}
Then, by the constructions of $u$ and $\tilde u$, one has
\begin{equation}\label{fracvetv}
  v(t)=\mathfrak{v}(t)\quad \text{for every } t\in [0,r_1].
\end{equation}
Then,  using \eqref{Delay1} together with the definition of
$T_{\gamma_1 \cdots \gamma_k}$ and $V(T_{\gamma_1 \cdots \gamma_k})$, one has
\begin{gather}\label{vandVTgamma}
v(T_{\gamma_1  \cdots \gamma_k} )=V(T_{\gamma_1 \cdots \gamma_k})
\quad \text{if }
T_{\gamma_1 \cdots \gamma_k}\in [0,T],
\end{gather}
%??? donner d\'{e}tails
with the convention that, if $k=0$, $T_{\gamma_1 \cdots \gamma_k}=T$.

Differentiating \eqref{Delay1} with respect to $t$, we get
\begin{equation}
\label{equation-der-1}
\big(1 + v_2'(t) \big)  v'\Big(t + r_2 + v_2(t) \Big)  = \big(1 + v_2'(t) \big) v_1'\Big(t + r_2 + v_2(t) - r_1\Big) G_1 + v_2'(t) G_2.
\end{equation}
It follows that
\begin{equation}
\label{expressionv'}
v'\Big(t + r_2 + v_2(t) \Big)  =  v_1'\Big(t + r_2 + v_2(t) - r_1\Big) G_1 + v_2'(t) G_2 - \frac{v_2'(t)^2}{1 + v_2'(t)} G_2.
\end{equation}
 From the definition of $dV$, \eqref{def-v-1}, \eqref{dVpetit},  \eqref{fracvetv}  and \eqref{expressionv'}, one gets, for every $T>r_1$,
 the existence of $C(T)>0$ such that
\begin{equation}\label{diffv'dV}
  | v'(T) - dV(T)| \le C(T) \eps^2.
\end{equation}
provided that $\eps$ is small enough (the smallness depending on $T$). In \eqref{diffv'dV} and in the following we use the notation
\begin{equation}\label{def|x|}
  |x|:=\|x \|_2 \quad \forall  x\in \R^n.
\end{equation}
 From \eqref{system},
\eqref{udecroitexp} and \eqref{deffrakv},
\begin{equation}
\label{expdecrea}
|v'(t)| \le 2\Lambda_2 C_0 e^{-\nu t} \|u^0 \|_{C^1([0,1];\R^2)}\quad \text{for every } t\in [0,+\infty),
\end{equation}
provided that $\|u^0 \|_{C^1([0,1];\R^2)}\leq \eps_0$. Using  \eqref{initVetdV}, \eqref{choice-u0}, \eqref{diffv'dV} and \eqref{expdecrea}, one gets  the existence of $C_1>0$ such that, for every $T>0$, there exist $C(T)>0$ and $\varepsilon(T)>0$ such that
\begin{equation}\label{estimate-contradiction}
  1\leq C_1 e^{-\nu T} + C(T) \eps \quad \text{for every } T>0,\text{ for every } \epsilon \in (0,\varepsilon(T)] .
\end{equation}
We choose $T>0$ large enough so that $C_1 e^{-\nu T}\leq (1/2)$. Then letting $\eps \rightarrow 0^+$ in \eqref{estimate-contradiction} we get a contradiction.

\medskip
It remains to prove \eqref{T} in order to conclude the proof of Theorem~\ref{thm1} if $m=1$. Let us assume
\begin{equation}\label{T-T}
T_{\gamma_1  \cdots \gamma_{k}} = T_{\alpha_1  \cdots \alpha_{m}} \text{ with } k,m\in \{1,\ldots,n-1\}
\end{equation}
($\gamma_i, \alpha_i =1, 2$).
Using \eqref{independence} and \eqref{estTgamma},  we derive that
\begin{equation}\label{ell}
m=k, \,  \mbox{card}\big\{i; \gamma_i = 2 \big\} = \mbox{card} \big\{i; \alpha_i = 2 \big\} =: \ell
\end{equation}
for some $0 \le \ell \le m $.
%The conclusion of the claim is now a consequence of \eqref{claimP}, \eqref{def-t1}, and \eqref{def-t2}.
Let $k_1 <   \cdots <    k_\ell$  and  $m_1 <    \cdots <   m_\ell$ be such that
\begin{equation*}
\gamma_{k_l} = \alpha_{m_l} = 2 \quad \mbox{ for } 1\le l \le \ell.
\end{equation*}
Define
\begin{equation*}
i_l := \sum_{i=1}^{k_l} (\gamma_i - 1) \quad \mbox{ and } \quad j_l := \sum_{i=1}^{k_l} (\alpha_i - 1).
\end{equation*}
It follows from \eqref{def-t2}, \eqref{def-l}, and \eqref{T-T} that
\begin{equation}\label{T-T-1}
\sum_{l = 1}^\ell t_{i_l}^{n - k_l} = \sum_{l = 1}^\ell t_{j_l}^{n- m_l}.
\end{equation}
Hence
\begin{equation}\label{conclusion-T-T}
\gamma_i = \alpha_i \quad \mbox{ for } i =1, \cdots, k=m
\end{equation}
 is proved if  one can verify that
\begin{equation}\label{claimP}
i_l = j_l \quad \mbox{ and } \quad k_l = m_l \quad \quad \forall \, l =1, \cdots \ell.
\end{equation}
By a recurrence argument on $\ell$, it suffices to prove that
\begin{equation}\label{claim1}
i_\ell = j_\ell \quad \mbox{ and } k_\ell = m_\ell.
\end{equation}

\medskip
Note that, by \eqref{construction-st},
\begin{equation}\label{observation}
t_{j}^k = a^k \eta^k t_{j + k}^0 + P_{k-1}(\xi, \eta),
\end{equation}
where $P_{k-1}$ is a polynomial of degree $k-1$ with rational coefficients.
Since $\xi, \eta$ satisfy \eqref{cond-xieta},  it follows from \eqref{T-T-1} and \eqref{observation} that
\begin{equation*}
k_\ell = m_\ell,
\end{equation*}
and
\begin{equation*}
i_\ell = j_\ell.
\end{equation*}
Thus claim \eqref{claim1} is proved and so are claims \eqref{claimP}, \eqref{conclusion-T-T}, and \eqref{T}. This concludes the proof of Theorem~\ref{thm1} if $m=1$.

Let us show how to modify the above proof to treat the case $m\geq 2$. Instead of \eqref{initial}, one requires

\begin{equation}\label{initial-new}
\|(s_i^0, t_i^0)  \|_\infty \le \eps^{2+m}/ 4^n \quad  \mbox{ for every } i, \,j \in \{1,\cdots,n\}.
\end{equation}
Then, instead of \eqref{norm-v}, one gets
\begin{equation}\label{norm-v-new}
\| \mathfrak{v}(T_{\alpha_1 \cdots \alpha_k})\|_\infty \le \eps^{2+m} \quad \mbox{ if } T_{\alpha_1 \cdots \alpha_k} \in (0, r_1).
\end{equation}
Instead of \eqref{def-v-1}, one requires
\begin{equation}\label{def-v-1-new}
\mathfrak{v}^{(m)} (T_{\alpha_1 \cdots \alpha_k}) = dV(T_{\alpha_1 \cdots \alpha_k}),
\end{equation}
and instead of \eqref{v-1}, one has
\begin{gather}
\label{v-1-new}
\|\mathfrak{v} \|_{C^m([0, r_1])} \le C \max\{\eps^2, A \},
\end{gather}
where $A$ is still given by \eqref{defA}. Then \eqref{choice-u0} is now
\begin{equation}\label{choice-u0-new}
\| u^0 \|_{C^m([0,1];\R^2)} \le C \| v \|_{C^m([0, r_1];\R^2)} \le C  \eps c^{T/ (2 r_1)}.
\end{equation}
 In the case $m=1$ we differentiated once \eqref{Delay1} with respect to $t$
 in order to get \eqref{expressionv'}. Now we differentiate
 \eqref{Delay1} $m$ times with respect to $t$ in order to get

\begin{gather*}
%\label{expressionv(m)}
\left|v^{(m)}\Big(t + r_2 + v_2(t) \Big)  -v_1^{(m)}\Big(t + r_2 + v_2(t) - r_1\Big) G_1 + v_2^{(m)}(t) G_2
\right|\leq
C \sum_{i=0}^m v^{(i)}(t)^2,
\end{gather*}
which allows us to get, instead of \eqref{diffv'dV},
\begin{equation}\label{diffv'dV-new}
  | v^{(m)}(T) - dV(T)| \le C(T) \eps^2.
\end{equation}
%\begin{equation}\label{diffv'dV-new}
%  \| v^{(m)}(T) - dV(T)\|_\infty \le C  \max\big\{\eps^2,  \eps c^{T/ (2 r_1)} \big\}.
%\end{equation}
We then get a contradiction as in the case $m=1$. This concludes the proof of Theorem~\ref{thm1}.

\proofend
\begin{remark}
\label{remdiffTcapital}
Property \eqref{T} is a key point. It explains why the condition $\hat \rho _0(K)<1$ is not sufficient
for exponential stability in the case of \textbf{nonlinear} systems. Indeed $\hat \rho _0(K)<1$ gives an exponential stability
which is robust with respect to  perturbations on the delays which are \textbf{constant}: these perturbations are not
allowed to depend on time. However with these type of perturbations \eqref{T} does not hold: with constant perturbations on the delays, one has
\begin{equation*}
%\label{Tgammaegaux}
  T_{12}=T_{21}, \, T_{122}=T_{212}=T_{221}
\end{equation*}
and, more generally,
\begin{equation*}
  T_{\gamma_1  \cdots  \gamma_k}=T_{\alpha_1  \cdots   \alpha_k}
  \text{ if }\mbox{card}\{i\in \{1,\cdots,k\};\, \gamma _i=1\}=
 \mbox{card}\{i\in \{1,\cdots,k\};\, \alpha_i=1\}.
\end{equation*}
\end{remark}

\section{Proof of Theorem~\ref{thm2}}\label{sect-thm2}

This section containing two subsections is devoted to the proof of Theorem~\ref{thm2}. In the first subsection, we present some lemmas which will be used in the proof. In the second subsection, we give the proof of Theorem~\ref{thm2}.

\subsection{Some useful lemmas}

The first lemma is standard one on the well-posedness of \eqref{sys-P} and \eqref{bdry-P}.

\begin{lemma}\label{lem1} Let $p\in [1,+\infty]$. There exist $C>0$ and $\gamma >0$ such that, for every $T>0$,
there exists $\eps_0>0$ such that,
for every $u_0 \in W^{2, p}((0, 1);\R^n)$ with
$\| u_0\|_{W^{2, p}((0, 1);\R^n)} < \eps_0$
satisfying the compatibility conditions
\eqref{compatibilty-C1-0}-\eqref{compatibilty-C1-1},
there exists one and only one solution $u\in C^1([0,T]\times [0,1];\R^n)$ of \eqref{system} satisfying the initial condition $u(0,\cdot)=u^0$. Moreover
\begin{equation*}
\| u(t, \cdot) \|_{W^{2, p}((0, 1);\R^n)}  \le C e^{\gamma t} \| u^0 \|_{W^{2, p}((0, 1);\R^n)}.
\end{equation*}
% As a consequence, for any $\delta_0 > 0$, there exists $\eps_0 > 0$ independent of $u_0$ such that if $ \| u_0\|_{W^{2, p}} \le \eps_0$ then
%\begin{equation*}
%\| u\|_{L^\infty([0, T] \times [0, 1])} + \| \partial_t u\|_{L^\infty([0, T] \times [0, 1])} + \| \partial_{x} u\|_{L^\infty([0, T] \times [0, 1])} \le \delta_0.
%\end{equation*}
\end{lemma}

We next present  two lemmas dealing with the system
\begin{equation*}
v_t+ A(t, x) v_x = 0,
\end{equation*}
and its perturbation where $A$ is diagonal. The first
lemma is the following one.

\begin{lemma}\label{lemP1} Let $p\in [1,+\infty]$, $m$ be a
positive integer,  $\lambda_1 \ge \cdots \ge \lambda_m  > 0$
and $\hat K\in (0,1)$. Then there exist three constants
 $\varepsilon_0>0$, $\gamma >0$ and $C>0$ such that, for every
 $T>0$, every $A \in C^1([0, T] \times [0, 1];\mathcal{D}_{m,+} )$,
 every $K \in C^1([0, T];\mathcal{M}_{m,m}(\R))$,
 every $v \in W^{1, p}([0, T] \times [0, 1]; \R^m)$
 such that
\begin{gather}
\label{eqvlinear}
v_t + A(t, x) v_x = 0 \mbox{ for } (t, x) \in   (0, T) \times (0, 1),
\\
\label{boundaryvlinear}
v(t, 0)  = K(t) v(t, 1)\mbox{ for } t \in   [0, T],
\\
\label{pro-K}
\sup_{t \in [0, T]} \|K(t) \|_p \le \hat K < 1,
\\
\label{derAK}
\| A - \mbox{diag}(\lambda_1, \cdots, \lambda_m)\|_{C^1([0, T]
\times [0, 1];\mathcal{M}_{m,m}(\R))} + \mathop{\sup}_{t \in [0, T]}\|K'(t) \|_{p} \le \eps_0,
\end{gather}
one has
\begin{equation*}
\|v(t, \cdot)\|_{W^{1, p}((0, 1);\R^m)} \le C e^{-\gamma t} \|v(0, \cdot)\|_{W^{1, p}((0, 1);\R^m)}
\mbox{ for } t \in   [0, T].
\end{equation*}
\end{lemma}

\noindent{\bf Proof of Lemma~\ref{lemP1}.}  We only consider the case $1 \le p < + \infty$,
the case $p=+ \infty$ follows similarly (the proof is even easier) and is left to the reader.
 For $t \ge 0$,
let $\varphi_i(t, s)$ be such that
\begin{equation*}
\partial_s \varphi_i(t, s) = A_{ii}(s, \varphi_i(t, s))  \quad \mbox{ and } \quad \varphi_i(t, t) = 0.
\end{equation*}
Then
\begin{equation*}
v_i(s, \varphi_i(t, s)) =v_i(t, 0).
\end{equation*}
We define  $s_i$ as a function of $t$  by   $\varphi_i (t, s_i(t)) = 1$.
Note that $A_{ii}(s, \varphi_i(t, s))>
\lambda_m/ 2 > 0$, at least if $\varepsilon_0>0$ is small enough, a property which is always assumed in this proof.
Hence $s_i$ is well-defined. It follows from the definition of $s_i$ that
\begin{equation}\label{def-si-lem}
v_i(s_i(t), 1) = v_i(t,0).
\end{equation}
%We have
%\begin{equation*}
%0 = \partial_t \varphi_i (t, s_i(t)) + \partial_s \varphi_i (t, s_i(t)) s_i'(t) =  \partial_t \varphi_i (t, s_i(t)) + A_{ii}(s, \varphi_i(t, s)) s_i'(t).
%\end{equation*}
%Here and in what follows $'$ denotes the derivative with respect to $t$, e.g., $s_i'(t) = ds_i/dt$ and $u'(t, x) = \partial_t v(t, x)$.
%On the other hand,
%\begin{equation*}
%\partial_s \partial_t \varphi_i(t, s) = \partial_x A_{ii}(s, \varphi_i(t, s)) \partial_t \varphi_i(t, s)
%\end{equation*}
%and
%\begin{equation*}
%\partial_t \varphi_i (t, t) =0.
%\end{equation*}
%Since $|\partial_x A_{ii}| \le \eps_0$, it follows that
%\begin{equation*}
%|\partial_t \varphi_i(t, s_i(t)) = \int_t^{s_i(t)} \partial_x A_ii(s, \varphi_i)
%\end{equation*}
Using classical results on the dependence of solutions of ordinary differential equations
on the initial conditions together with the inverse mapping theorem, one gets
\begin{equation}\label{dsi-1-lem}
|s_i'(t) - 1| \le C \eps_0.
\end{equation}
Here and in what follows in this proof $'$ denotes the derivative with
respect to $t$, e.g., $s_i'(t) = ds_i/dt$ and $v'(t, x) = \partial_t v(t, x)$ and
$C$ denotes a positive constant which changes from one
place to another and may depend on $p$, $m$, $\lambda_1 \ge \cdots \ge \lambda_m  > 0$
and $\hat K\in (0,1)$ but is independent of $\eps_0>0$, which is always
assumed to be small enough, $T>0$, $A$ and $v$ which are always assumed to
satisfy \eqref{eqvlinear} to \eqref{derAK}.

\medskip
%Since $\|A - \mbox{diag} (\lambda_1, \cdots, \lambda_m)\|_{W^{2, p} \big((0, T) \times (0, 1) \big)}$ is small, it follows  that  $s_i$ is invertible (strictly increasing).
Define, for $t \ge 2 \lambda_1$,
\begin{equation}\label{def-ri-lem}
\hat r_i(t) :=  t- s_i^{-1}(t).
\end{equation}
%\begin{equation}\label{def-ri}
%\tilde r_i(t) = F_{ii} \Big( u \big(s_i^{-1}(t), 0 \big) \Big)
%\end{equation}
From \eqref{dsi-1-lem}, we have
%\begin{equation}\label{derivative-ri1-lem}
%(s_i^{-1})'(\xi) = A_{ii}(s_i(t), 1)/ A_{ii}(t, 0),
%\end{equation}
%where $t = s_i^{-1}(\xi)$. It follows from \eqref{derivative-ri1-lem} that
\begin{equation}\label{derivative-ri-11-lem}
\sup_{t \in [2 \lambda_1, T]} | \hat r_i'|\le C \eps_0.
\end{equation}
%and a combination of \eqref{assumption-lem} and \eqref{derivative-ri2-lem} implies
%\begin{equation}\label{derivative-ri-21-lem}
%\int_{2r_1}^{T} |\tilde r_i''(t)|_p^p \le C \int_{0}^{T} \|u''(t)\|_p^p + C \delta_1^p  \int_{0}^{T} \|u'(t) \|_p^p
%\end{equation}
Set
\begin{equation*}
V(t) = v(t, 0).
\end{equation*}
We derive from \eqref{boundaryvlinear}, \eqref{def-si-lem} and \eqref{def-ri-lem} that
\begin{equation}\label{important-lem}
V(t) = K(t) \Big(V_1\big(t - \hat r_1 (t)\big), \cdots, V_i \big(t - \hat r_i (t)\big), \cdots, V_m\big(t - \hat r_m (t)\big)\Big)^T, \quad \mbox{ for } t \ge 2 r_m.
\end{equation}
In \eqref{important-lem}  and in the following $r_i:=1/\lambda_i$ for every $i\in \{1,\cdots,m\}$.
 From  \eqref{pro-K} and \eqref{important-lem}, we obtain
\begin{equation}\label{key-1-lem}
\int_{2 r_m}^{T}\| V(t) \|_{p}^{p} \, dt \le \hat K^p \sum_{i=1}^{n} \int_{2 r_m}^T |V_{i} \big(t - \hat r_{i}(t) \big)|^{p} \, dt.
\end{equation}
Since
\begin{equation*}
\int_{2 r_m}^T |V_{i} \big(t - \hat r_{i} (t) \big)|^{p} \, dt =  \int_{2 r_m - \hat r_i (2 r_m) }^{T - \hat  \lambda_i(T)} |V_{i}(t)|^{p} s_i' (t) \, dt,
\end{equation*}
it follows from \eqref{dsi-1-lem} that
\begin{equation}\label{inter-1-lem}
\int_{2 r_m}^T |V_{i}(t - \hat r_{i})|^{p} \le \int_{0}^{T} (1 + C\eps_0) |V_{i}(t)|^{p} \, dt.
\end{equation}
A combination of \eqref{key-1-lem} and \eqref{inter-1-lem} yields
\begin{equation*}
\int_{2 r_m}^T \| V(t) \|_p^p \, dt \le \int_0^{T} \hat K^p (1 + C \eps_0)\|V(t) \|_p^p \, dt.
\end{equation*}
By taking  $\eps_0$ small enough so that  $\hat K^p (1 + C \eps_0) \le [(1 + \hat K)/2 ]^p$, we have
\begin{equation}\label{fact1-lem}
\int_0^T \|V(t) \|_p^p \, dt \le C \int_0^{2 r_m} \|V(t) \|_p^p \, dt.
\end{equation}

We next establish  similar estimates for the derivatives of $V$. Let us define
\begin{equation}\label{def-w-lem}
W(t) : =(W_1(t),\cdots,W_m(t))^T:=V'(t).
\end{equation}
Differentiating \eqref{important-lem} with respect to  $t$, we have
\begin{equation}\label{eq-w-lem}
W(t) =   K(t) \Big(W_1\big(t - \hat r_1 (t)\big), \cdots, W_i \big(t - \hat r_i (t)\big) , \cdots, W_m\big(t - \hat r_m (t)\big) \Big)^T  + g_1(t) + f_1(t),
\end{equation}
where
\begin{equation}\label{def-g1-lem}
g_1(t) := - K(t) \Big(W_1\big(t - \hat r_1 (t)\big) \hat r_1'(t) , \cdots, W_i \big(t - \hat r_i (t)\big) \hat r_i'(t) , \cdots, W_m\big(t - \hat r_m (t)\big) \hat r_m'(t) \Big)^T
\end{equation}
and
\begin{equation}\label{def-f1-lem}
f_1(t) := K'(t) \Big(V_1\big(t - \hat r_1 (t)\big), \cdots, V_i \big(t - \hat r_i (t)\big), \cdots, V_m\big(t - \hat r_m (t)\big)\Big)^T.
\end{equation}
From \eqref{eq-w-lem}, we have
\begin{equation}\label{estimate-dv}
|W(t)|_p^p \le [(\hat K + 1)/2]^p \sum_{i=1}^{m} |W_{i} \big(t - \hat r_{i}(t) \big)|^{p} + C\Big(|f_1(t) |_p^p +  |g_1(t)|_p^p \Big).
\end{equation}
Using \eqref{derAK} and \eqref{derivative-ri-11-lem}, we derive from \eqref{def-g1-lem} and \eqref{def-f1-lem}, as in \eqref{inter-1-lem}, that
\begin{equation}\label{estimate-g1-f1-lem}
\int_{2 r_m}^T \big(\|g_1(t) \|_p^p  + \|f_1(t) \|_p^p\big) \, dt \le C\eps_0^p \int_{0}^{T} \big( \|W \|_p^p + \|V(t) \|_p^p \big)\, dt.
\end{equation}
It follows from \eqref{estimate-dv}, as in \eqref{fact1-lem},  that
\begin{equation}\label{fact2-lem}
\int_0^T \|V'(t) \|_p^p \, dt \le C \int_0^{2 r_m} \big( \|V(t) \|_{p}^p + \|V'(t) \|_p^p \big)\, dt.
\end{equation}
Combining \eqref{fact1-lem} and \eqref{fact2-lem}, we reach the conclusion.  \proofend

\medskip
As a consequence of Lemma~\ref{lemP1}, we obtain the following lemma, where
$\mathcal{B}(\R^m)$ denotes the set of bilinear forms on $\R^m$.
\begin{lemma}\label{lemP2} Let $p \ge 1$, $m$ be a
positive integer,  $\lambda_1 \ge \cdots \ge \lambda_m  > 0$,
$\hat K\in (0,1)$ and $M\in(0,+\infty)$. Then there exist three constants
 $\varepsilon_0>0$, $\gamma >0$ and $C>0$ such that, for every
 $T>0$, every $A \in C^1([0, T] \times [0, 1];\mathcal{D}_{m,+} )$,
 every $K \in C^1([0, T];\mathcal{M}_{m,m}(\R))$, every
 $Q \in C^1([0, T] \times [0, 1];\mathcal{B}(\R^m) )$
  and every $v \in W^{1, p}([0, T] \times [0, 1]; \R^m)$
 such that
\begin{gather}
\label{vequation}
v_t + A(t, x) v_x = Q(t, x)(v, v) \mbox{ for } (t, x) \in   (0, T) \times (0, 1),
\\
\label{vbord}
v(t, 0)  = K(t) v(t, 1)\mbox{ for } t \in   (0, T),
\\
\label{Kp<1}
\sup_{t \in [0, T]}\|K(t) \|_p \le \hat K < 1,
\\
\| A - \mbox{diag}(\lambda_1, \cdots, \lambda_m)\|_{C^1([0, T] \times [0, 1])}
+ \mathop{\sup}_{t \in [0, T]}\|K'(t) \|_{p}\le \eps_0,
\\
\|Q\|_{C^1([0, T] \times [0, 1];\mathcal{B}(\R^m) )}\leq M,
\\
\label{V0petit}
\|v(0, \cdot)\|_{W^{1, p}((0, 1);\R^m)}\le \eps_0,
\end{gather}
one has
\begin{equation*}
\|v(t, \cdot)\|_{W^{1, p}((0, 1);\R^m)} \le C e^{-\gamma t}
 \|v(0, \cdot)\|_{W^{1, p}((0, 1);\R^m)} \mbox{ for t} \in (0,T).
\end{equation*}
\end{lemma}
\noindent{\bf Proof of Lemma~\ref{lemP2}.}
Let
$\tilde v \in W^{1, p}([0, T] \times [0, 1]; \R^m)$ be the solution of the linear
Cauchy problem
\begin{gather}\label{deftildev-eq}
  \tilde v_t + A(t, x) \tilde v_x = 0
  \mbox{ for } (t, x) \in   (0, T) \times (0, 1),
\\
\label{deftildev-bord}
\tilde v(t, 0)  = K(t) \tilde v(t, 1)\mbox{ for } t \in   (0, T),
\\
\label{deftildev-init}
\tilde v(0, x)  = v(0, x) \mbox{ for } x \in   (0, 1).
\end{gather}
(Note that $v(0,0)=K(0)v(0,1)$; hence such a $\tilde v$ exists.)  From Lemma~\ref{lemP1},
\eqref{deftildev-eq},
 \eqref{deftildev-bord} and \eqref{deftildev-init}, one has
\begin{equation}\label{estsurtildev}
\|\tilde v(t, \cdot)\|_{W^{1, p}((0, 1);\R^m)}
\le C e^{-\gamma t} \|v(0, \cdot)\|_{W^{1, p}((0, 1);\R^m)}
\mbox{ for } t \in   [0, T].
\end{equation}
Let
\begin{gather}\label{defbarv}
\bar v:=v-\tilde v.
\end{gather}
 From \eqref{vequation}, \eqref{vbord}, \eqref{deftildev-eq},
 \eqref{deftildev-bord}, \eqref{deftildev-init} and \eqref{defbarv}, one has
\begin{gather}\label{barv-eq}
 \bar v_t + A(t, x) \bar v_x = Q(t, x)(\tilde v+\bar v , \tilde v+ \bar v)
  \mbox{ for } (t, x) \in   (0, T) \times (0, 1),
\\
\label{barv-bord}
\bar v(t, 0)  = K(t)  \bar v(t, 1)\mbox{ for } t \in   (0, T),
\\
\label{barv-init}
\bar v(0, x)  = 0 \mbox{ for } x \in   (0, 1).
\end{gather}
Let, for $t\in [0,T]$,
\begin{equation}\label{defe(t)}
  e(t):=\|\bar v(t,\cdot)\|_{L^\infty((0,1);\R^m)}.
\end{equation}
Following the characteristics and using \eqref{estsurtildev}, \eqref{barv-eq}, \eqref{barv-bord} and
the Sobolev imbedding $W^{1, p}((0, 1);\R^m)\subset L^\infty((0,1);\R^m)$, one gets, in the
sense of distribution in $(0,T)$,
\begin{equation}\label{este(t)}
  e'(t)\leqslant C (\|v(0, \cdot)\|_{W^{1, p}((0, 1);\R^m)}^2+ e(t)+ e(t)^2).
\end{equation}
In \eqref{este(t)},  $C$ is as in the proof of
Lemma~\ref{lemP1} except that it may now depend on $M$. From \eqref{barv-init}, \eqref{defe(t)} and \eqref{este(t)}, one gets
 the existence of $\varepsilon_0$, of an increasing function $T\in[0,
 +\infty) \mapsto C(T)\in (0,+\infty)$ and of a decreasing function $T\in[0,
 +\infty) \mapsto \varepsilon(T)\in (0,+\infty)$, such that, for every
$T\in [0,+\infty)$, for every
$A \in C^1([0, T] \times [0, 1];\mathcal{D}_{m,+} )$,
 every $K \in C^1([0, T];\mathcal{M}_{m,m}(\R))$, every
 $Q \in C^1([0, T] \times [0, 1];\mathcal{B}(\R^m) )$
  and every $v \in W^{1, p}([0, T] \times [0, 1]; \R^m)$ satisfying \eqref{vequation}
  to \eqref{V0petit},
\begin{multline}\label{epetit}
  \left(\|v(0, \cdot)\|_{W^{1, p}((0, 1);\R^m)}\le \eps(T)\right)\implies
\\
\left( \|\bar v(t,\cdot)\|_{L^\infty((0,1);\R^m)} \le C (T)
 \|v(0, \cdot)\|_{W^{1, p}((0, 1);\R^m)}^{2} \mbox{ for } t \in (0,T)\right),
\end{multline}

Let $\bar w:=\bar v_x$.
Differentiating
 \eqref{barv-eq} with respect to $x$, we get
\begin{multline}\label{barv-eq-der}
 \bar w_t + A(t, x) \bar w_x + A_x(t,x) \bar w= Q_x(t, x)(\tilde v+\bar v , \tilde v+ \bar v)
 \\+Q(t, x)(\tilde v_x+\bar w , \tilde v+ \bar v)
 + Q(t, x)(\tilde v+\bar v , \tilde v_x+ \bar w)
  \mbox{ for } (t, x) \in   (0, T) \times (0, 1).
\end{multline}
Differentiating \eqref{barv-bord} with respect to $t$ and using \eqref{barv-eq}, we get,
for  $t \in   [0, T]$,
\begin{multline}\label{barv-bord-der}
A(t, 0)\bar w(t,0) - Q(t, 0)(\tilde v(t,0)+\bar v (t,0), \tilde v(t,0)+ \bar v(t,0))
  =
\\K(t) \big(A(t, 1)\bar w(t,1) - Q(t,1)(\tilde v(t,1)+\bar v (t,1), \tilde v(t,1)+ \bar v(t,1)) \big)
  -K'(t)\bar v(t,1).
\end{multline}
Differentiating \eqref{barv-init} with respect to $x$, one gets
\begin{gather}
\label{barw-init}
\bar w(0, x)  = 0 \mbox{ for } x \in   (0, 1).
\end{gather}
We consider \eqref{barv-eq-der}, \eqref{barv-bord-der} and \eqref{barw-init} as a nonhomogeneous linear hyperbolic
system where the unknown is $w$ and the data are $A$, $K$, $Q$, $\tilde v $,  and $\bar v$.  Then,
from straightforward estimates on the solutions of linear hyperbolic equations, one gets that, for every $t\in [0,T]$,
\begin{equation}\label{estimatewLp}
\begin{array}{rcl}
\displaystyle
\|\bar w(t,\cdot)\|_{L^p((0,1);\R^m)}&\leq &e^{CT\left(1+ \|\tilde v\|_{L^\infty((0,T)\times(0,1);\R^m)}+
\|\bar v\|_{L^\infty((0,T)\times(0,1);\R^m)}\right)}
\\
&&
\displaystyle\times\left(\|\tilde v\|_{L^\infty((0,T);W^{1,p}((0,1);\R^m))}^2+\|\bar v\|_{L^\infty((0,T)\times(0,1);\R^m)}^2\right).
\end{array}
\end{equation}
From \eqref{estsurtildev}, \eqref{epetit} and  \eqref{estimatewLp}, one gets
 the existence of $\varepsilon_0$, of an increasing function $T\in[0,
 +\infty) \mapsto C(T)\in (0,+\infty)$ and of a decreasing function $T\in[0,
 +\infty) \mapsto \varepsilon(T)\in (0,+\infty)$, such that, for every
$T\in [0,+\infty)$, every
$A \in C^1([0, T] \times [0, 1];\mathcal{D}_{m,+} )$,
 every $K \in C^1([0, T];\mathcal{M}_{m,m}(\R))$, every
 $Q \in C^1([0, T] \times [0, 1];\mathcal{B}(\R^m) )$
  and every $v \in W^{1, p}([0, T] \times [0, 1]; \R^m)$ satisfying \eqref{vequation}
  to \eqref{V0petit},
\begin{multline}\label{v0petit0K}
  \left(\|v(0, \cdot)\|_{W^{1, p}((0, 1);\R^m)}\le \eps(T)\right)\implies
\\
\left(\|\bar v(t, \cdot)\|_{W^{1, p}((0, 1);\R^m)} \le C (T)
 \|v(0, \cdot)\|_{W^{1, p}((0, 1);\R^m)}^{2} \mbox{ for } t \in (0,T)\right),
\end{multline}
which, together with \eqref{estsurtildev} and \eqref{defbarv},
concludes the proof of Lemma~\ref{lemP2}.
\endproof

\subsection{Proof of Theorem~\ref{thm2}}

Replacing, if necessary, $u$ by $Du$ where $D$ (depending only on $K$)
is a diagonal matrix with positive entries, we may
assume that
\begin{equation}\label{p-normG'0}
  \|G'(0)\|_p<1.
\end{equation}

%We assume that the system is in the characteristic form.
For $a\in \R^n$, let $\lambda_i(a)$ be the $i$-th eigenvalue of $F(a)$
and $l_i(a)$ be a left eigenvector of $F(a)$ for this eigenvalue.
The functions $\lambda_i$ are of class
$C^\infty$ in a neighborhood of
$0\in \R^n$. We may also impose on the $l_i$ to
be of class $C^\infty$ in a neighborhood
of $0\in \R^n$ and that $l_i(0)^T$ is the $i$-th vector
of the canonical basis of $\R^n$.
Set
\begin{equation*}
\left\{\begin{array}{l}
v_i = l_i(u) u \\[6pt]
w_i = l_i(u) \partial_t u
\end{array} \right.  \quad \mbox{ for } i =1, \cdots, n.
\end{equation*}
From \cite[(3.5) and (3.6) on page 187]{Li-book}, we have, for  $i =1, \cdots, n$,
\begin{equation}\label{sys1}
\left\{\begin{array}{l}
u_i = v_i + \sum_{j, k}^n b_{ijk}(v) v_j v_k\\[6pt]
\partial_t u_i = w_i + \sum_{ijk} \bar b_{ijk} (v) v_j w_k
\end{array} \right. ,
\end{equation}
where $b_{ijk}$ and  $\bar b_{ijk}$ are of class $C^\infty$.
From \cite[(3.7) and (3.8)]{Li-book}, we obtain,
for  $i =1, \cdots, n$,
\begin{equation}\label{sys2}
\left\{\begin{array}{l}
\dsp \partial_t v_i + \lambda_i(u) \partial_x v_i
= \sum_{ijk}^n c_{ijk}(u)v_j v_k + \sum_{ijk}^n d_{ijk}(u) v_j w_k, \\[6pt]
\dsp \partial_t w_i + \lambda_i(u) \partial_x w_i = \sum_{ijk}^n \bar c_{ijk}(u)w_j w_k + \sum_{ijk}^n \bar d_{ijk}(u) v_j w_k, \\[6pt]
\end{array}\right.
\end{equation}
where $c_{ijk}, \bar c_{ijk}, d_{ijk}, \bar d_{ijk}$ are of
class $C^\infty$ in a neighborhood of $0\in \R^n$. We also have,
for some $\hat G: \mR^{2n} \to \mR^{2n}$ of class $C^\infty$ in a neighborhood of $0\in \R^{2n}$,
\begin{equation*}
\left( \begin{array}{c} v(t, 0) \\[6pt]
w(t, 0)
\end{array}
\right) = \hat G \left( \begin{array}{c} v(t, 1) \\[6pt]
w(t, 1)
\end{array}
\right)
\end{equation*}
and, by \eqref{bdry-P},
\begin{equation*}
\hat G'  \left( \begin{array}{c} 0 \\[6pt]
0
\end{array}
\right) =  \left( \begin{array}{cc }  G'(0) & 0 \\[6pt]
0 & G'(0)
\end{array}
\right),
\end{equation*}
which, together with \eqref{p-normG'0}, implies that
\begin{equation*}
\|\hat G'(0)\|_p < 1.
\end{equation*}
Applying  Lemma~\ref{lemP2} for \eqref{sys2}, we obtain the
exponential stability  for $(v, w)$ with respect to the $W^{1, p}$-norm, from which, noticing that $u_x=-F(u)^{-1}u_t$,
Theorem~\ref{thm2} readily follows.
\proofend

\bibliographystyle{plain}

\bibliography{biblio-CE}

\def\cprime{$'$} \newcommand{\SortNoop}[1]{}
\begin{thebibliography}{10}

\bibitem{2014-Coron-Bastin}
Jean-Michel Coron and Georges Bastin.
\newblock Dissipative boundary conditions for one-dimensional quasi-linear
  hyperbolic systems: Lyapunov stability for the {$C^1$}-norm.
\newblock {\em Preprint}, 2014.

\bibitem{2008-Coron-Bastin-Novel-SICON}
Jean-Michel Coron, Georges Bastin, and Brigitte d'Andr{\'e}a Novel.
\newblock Dissipative boundary conditions for one-dimensional nonlinear
  hyperbolic systems.
\newblock {\em SIAM J. Control Optim.}, 47(3):1460--1498, 2008.

\bibitem{2007-Coron-Andrea-Novel-Bastin-IEEE}
Jean-Michel Coron, Brigitte d'Andr{\'e}a Novel, and Georges Bastin.
\newblock A strict {L}yapunov function for boundary control of hyperbolic
  systems of conservation laws.
\newblock {\em IEEE Trans. Automat. Control}, 52(1):2--11, 2007.

\bibitem{GreenbergLi}
James~M. Greenberg and Ta-tsien Li.
\newblock The effect of boundary damping for the quasilinear wave equation.
\newblock {\em J. Differential Equations}, 52(1):66--75, 1984.

\bibitem{HaleVerduynLunel-book}
Jack~K. Hale and Sjoerd~M. Verduyn~Lunel.
\newblock {\em Introduction to functional-differential equations}, volume~99 of
  {\em Applied Mathematical Sciences}.
\newblock Springer-Verlag, New York, 1993.

\bibitem{2003-De-Halleux-et-al-Automatica}
Jonathan {\SortNoop{Halleux}}de~Halleux, Christophe Prieur, Jean-Michel Coron,
  Brigitte {\SortNoop{Andr\'{e}a-Novel, Brigitte}}~{d'Andr\'{e}a-Novel}, and
  Georges Bastin.
\newblock Boundary feedback control in networks of open channels.
\newblock {\em Automatica J. IFAC}, 39(8):1365--1376, 2003.

\bibitem{Li-book}
Ta-tsien Li.
\newblock {\em Global classical solutions for quasilinear hyperbolic systems},
  volume~32 of {\em RAM: Research in Applied Mathematics}.
\newblock Masson, Paris, 1994.

\bibitem{LiYu}
Ta-tsien Li and Wen~Ci Yu.
\newblock {\em Boundary value problems for quasilinear hyperbolic systems}.
\newblock Duke University Mathematics Series, V. Duke University Mathematics
  Department, Durham, NC, 1985.

\bibitem{2008-Prieur-Winkin-Bastin-MCSS}
Christophe Prieur, Joseph Winkin, and Georges Bastin.
\newblock Robust boundary control of systems of conservation laws.
\newblock {\em Mathematics of Control, Signal and Systems (MCSS)}, 20:173--197,
  2008.

\bibitem{1985-Qin-Tie-hu}
Tie~Hu Qin.
\newblock Global smooth solutions of dissipative boundary value problems for
  first order quasilinear hyperbolic systems.
\newblock {\em Chinese Ann. Math. Ser. B}, 6(3):289--298, 1985.
\newblock A Chinese summary appears in Chinese Ann.\ Math.\ Ser.\ A {\bf 6}
  (1985), no.\ 4, 514.

\bibitem{1974-Rauch-Taylor-IUMJ}
Jeffrey Rauch and Michael Taylor.
\newblock Exponential decay of solutions to hyperbolic equations in bounded
  domains.
\newblock {\em Indiana Univ. Math. J.}, 24:79--86, 1974.

\bibitem{Slemrod}
Marshall Slemrod.
\newblock Boundary feedback stabilization for a quasilinear wave equation.
\newblock In {\em Control theory for distributed parameter systems and
  applications ({V}orau, 1982)}, volume~54 of {\em Lecture Notes in Control and
  Inform. Sci.}, pages 221--237. Springer, Berlin, 1983.

\bibitem{02Xu-Sallet}
Cheng-Zhong Xu and Gauthier Sallet.
\newblock Exponential stability and transfer functions of processes governed by
  symmetric hyperbolic systems.
\newblock {\em ESAIM Control Optim. Calc. Var.}, 7:421--442 (electronic), 2002.

\bibitem{1986-Zhao-Yan-chun}
Yan~Chun Zhao.
\newblock Classical solutions for quasilinear hyperbolic systems.
\newblock {\em Thesis, Fudan University}, 1986.
\newblock In Chinese.

\end{thebibliography}
%\bibliography{/Users/hnguyen/Dropbox/rho-optimial/V1/biblio-CE}

%\bibliographystyle{amsplain}
%\bibliography{/Dropbox/Bib/bib1}
%\bibliography{/Users/hoaiminhnguyen/Dropbox/Bib/bib1}
%\bibliography{/Users/hnguyen/Dropbox/Bib/bib1}

\end{document}